\documentstyle[amscd,amssymb,verbatim,diagrams]{amsart}
\newcommand{\vareps}{\varepsilon}

\newcommand{\Id}{\operatorname{Id}}

\newcommand{\can}{\operatorname{can}}

\renewcommand{\mod}{\operatorname{mod}}

\newcommand{\Cone}{\operatorname{Cone}}
\newcommand{\Tr}{\operatorname{Tr}}

\newcommand{\und}{\underline}

\newcommand{\OO}{{\cal O}}

\newcommand{\DD}{{\cal D}}

\newcommand{\BB}{{\cal B}}

\newcommand{\G}{{\Bbb G}}

\newcommand{\hra}{\hookrightarrow}

\newcommand{\lan}{\langle}
\newcommand{\ran}{\rangle}

\newcommand{\de}{\delta}
\newcommand{\eps}{\epsilon}

\renewcommand{\ker}{\operatorname{ker}}

\newcommand{\mf}{\operatorname{mf}}
\numberwithin{equation}{section}

\newtheorem{thm}{Theorem}[section]
\newtheorem{prop}[thm]{Proposition}
\newtheorem{lem}[thm]{Lemma}

\newenvironment{rem}{\vspace{3mm}\noindent
{\bf Remark.}}{\vspace{3mm}}

\newenvironment{ex}{\vspace{3mm}\noindent
{\bf Example.}}{\vspace{3mm}}

\newcommand{\Pf}{\noindent {\it Proof}}
\newcommand{\id}{\operatorname{id}}

\newcommand{\we}{\wedge}

\newcommand{\univ}{\operatorname{univ}}

\renewcommand{\AA}{{\cal A}}

\newcommand{\Om}{\Omega}
\newcommand{\dbar}{\overline{\partial}}
\newcommand{\Hom}{\operatorname{Hom}}

\newcommand{\Ext}{\operatorname{Ext}}
\newcommand{\End}{\operatorname{End}}

\renewcommand{\a}{\alpha}
\renewcommand{\b}{\beta}

\newcommand{\De}{\Delta}
\newcommand{\la}{\lambda}

\renewcommand{\H}{{\Bbb H}}

\newcommand{\Z}{{\Bbb Z}}
\newcommand{\Q}{{\Bbb Q}}

\newcommand{\Ga}{\Gamma}

\newcommand{\wt}{\widetilde}
\newcommand{\ot}{\otimes}

\newcommand{\sub}{\subset}
\newcommand{\ed}{\qed\vspace{3mm}}

\newcommand{\fg}{{\frak g}}

\newcommand{\tr}{\operatorname{tr}}

\newcommand{\unit}{{\mathbf 1}}

\newcommand{\ev}{\operatorname{ev}}

\newcommand{\at}{\operatorname{at}}
\newcommand{\MF}{\operatorname{MF}}
\newcommand{\str}{\operatorname{str}}

\title{Atiyah class and Chern character for global matrix factorizations}
\author{Bumsig Kim} 
\address{Korea Institute for Advanced Study}
\author{Alexander Polishchuk}
\address{University of Oregon, Korea Institute for Advanced Study
and National Research University Higher School of Economics, Russian Federation} 
\begin{document}

\begin{abstract}
We define the Atiyah class for global matrix factorizations and use it to give a formula for the categorical Chern character
and the boundary-bulk map for matrix factorizations, 
generalizing the formula in the local case obtained in \cite{PV-Chern}. Our approach
is based on developing the Lie algebra analogies observed by Kapranov \cite{Kapranov} and Markarian 
\cite{Markarian}.
\end{abstract}

\maketitle

\section*{Introduction}

Recall that for a vector bundle $E$ over a smooth algebraic variety $X$, one has a natural
class $\at(E)\in \Ext^1(E,\Om^1_X\ot E)$, called {\it Atiyah class} of $E$, such that
the Chern character of $E$ is obtained as $\tr(\exp(\at(E)))$ when $X$ is a projective smooth variety over $\mathbb{C}$ (see \cite{Atiyah}).
This construction also generalizes to bounded complexes of vector bundles (see \cite{Huybrechts Lehn}).
Furthermore, it shows up in the formula for the categorical Chern character and more generally, the boundary-bulk
map (see \cite{Caldararu-HKR}, \cite{Ram}).

The goal of this paper is to generalize the construction of the Atiyah class 
to the case of (global) matrix factorizations and to give a formula for the categorical boundary-bulk map for matrix factorizations. 

Thus, we start with a smooth scheme $X$ over $\mathbb{C}$ equipped with a function
$w$. We refer to \cite{EP}, \cite{LP}, \cite{Pre} for the background on categories of matrix factorizations.
 
We want to construct an analog of Atiyah class for a matrix factorization $(E,\de)$ of $w$.
In the case when $X$ is affine the construction of such an Atiyah class is known (see \cite{DM}, \cite{Yu}).

It turns out that in the case when $w\neq 0$ one can still construct a certain class $\hat{\at}(E)$ (see below)
which reduces to $1+\at(E)$ when $w=0$.

Below we denote by $\Hom^*_{\MF(w)}(\cdot,\cdot)$ the cohomology of the ($\Z/2$-graded) 
morphisms spaces in the category of matrix factorizations of $w$. When considering the Hochschild homology/cohomology
of $\MF(w)$,
we work over $k[u,u^{-1}]$, where $\deg(u)=2$ (using the $\Z/2$-grading on the category $\MF(w)$.

\medskip

\noindent
{\bf Main construction} (see Sec.\ \ref{global-atiyah-def-sec}). {\it 
For every matrix factorization $E$ of $w$ we construct a natural class
$$\hat{\at}(E)\in \Hom^0_{\MF(w)}(E, [\OO_X\rTo{dw}\Om^1_X]\ot E),$$
%\footnote{BK: in the derived category of $MF(X, w)$? I am here confused with the notation $\Hom ^0$ since it usually means the group of degree 0 maps in the dg category, which may not be closed maps.}
such that the image 
%\footnote{BK: Do you mean the composition?}
of  $\hat{\at}(E)$ under the natural projection 
$\Hom^0_{\MF(w)}(E,[\OO_X\rTo{dw}\Om^1_X]\ot E)\to \Hom^0_{\MF(w)}(E,E)$ is the identity element
$\id_E$. The formation of $\hat{\at}(E)$ is functorial and compatible with pull-backs.}

\bigskip

Next, we define the class 
$$\exp(\at(E))\in \Hom^0_{\MF(w)}(E,(\Om_X^\bullet,\we dw)\ot E)$$
as the composition of the iterated map 
$$\hat{\at}(E)^{(n)}:E\to [\OO_X\rTo{dw}\Om^1_X]^{\ot n}\ot E\to S^n[\OO_X\rTo{dw}\Om^1_X]\ot E,$$ 
where $n=\dim X$, with the isomorphism 
$$S^n[\OO_X\rTo{dw}\Om^1_X]\ot E\rTo{\sim} (\Om_X^\bullet,\we dw)\ot E$$
induced by the isomorphism
\begin{equation}\label{Sn-isom}
S^n[\OO_X\rTo{dw}\Om^1_X]=[\OO_X\rTo{n dw}\Om^1\rTo{(n-1)dw}\to\ldots]\to (\Om_X^\bullet,\we dw)
\end{equation}
given by $\a_0\mapsto \a_0$, $\a_i\mapsto \frac{\a_i}{n(n-1)\ldots (n-i+1)}$, where $\a_i\in\Om_X^i$.

The above definition may look a bit strange, however, it is easily explained by the fact that
when one tries to recover $\exp(x)$ from $(1+x)^n$ in the ring $\Q[x]/(x^{n+1})$, 
one has to rescale $x^i$ by the factor $\frac{1}{n(n-1)\ldots (n-i+1)}$. 

Below we view $\exp(\at(E))$ as an element of $\H^0$ of the $\Z/2$-graded complex $\und{\Hom}(E,E)\ot (\Om_X^\bullet,\we dw)$
and denote by
$$\str: \und{\Hom}(E,E)\ot (\Om_X^\bullet,\we dw)\to (\Om_X^\bullet,\we dw)$$
the supertrace morphism.
Here is our main result. 

\medskip

\noindent
{\bf Theorem A}. {\it Assume now that $w=0$ on the critical locus of $w$ (set-theoretically). 
Under the natural identification 
\begin{equation}\label{HKR-w-ThmA-eq}
HH_*(\MF(w))\simeq H^*(X,(\Om_X^\bullet,\we dw)),
\end{equation}
%\footnote{BK: Is $R\Gamma (X, (\Om_X^\bullet \we dw ))$ instead of $(\Om_X^\bullet \we dw )$?}
the categorical boundary-bulk map, 
%\footnote{BK: in $D(k)$?}
$$\Hom^*_{\MF(w)}(E,E)\to HH_*(\MF(w)),$$
for a matrix factorization $E$ of $w$, is equal to the map induced on hypercohomology by the map
$$\und{\Hom}(E,E)\to (\Om_X^\bullet,\we dw): x\mapsto \str(\exp(\at(E)) \cdot x )$$
in the $\Z/2$-graded derived category of $X$.
}

\medskip

Note that we give two constructions of isomorphism \eqref{HKR-w-ThmA-eq}: an abstract one coming from analogies with  Lie theory described below 
and the one given by an explicit chain map (this construction is due to \cite{LP}).

In the case of matrix factorizations over a regular local ring there is a simpler formula for the categorical Chern character
obtained in \cite{PV-Chern}. It can be derived from the above Theorem using a connection on the underlying vector bundle of
$E$ (see Remark \ref{trivial-bun-rem}). In the case of global matrix factorizations a formula similar to the one in Theorem A
was obtained by Platt \cite{Platt}. However, his definition of $\exp(\at(E))$ is much more complicated 
(based on some explicit resolutions of the relevant objects in the derived category).

Note that for $w=0$, i.e., in the classical case of vector bundles, we get a new proof of the compatibility of the categorical
Chern character with the classical one, originally proved by Caldararu \cite{Caldararu-HKR}. 
Our proof is more conceptual than that in
\cite{Caldararu-HKR}: we check the key techical statement (see Lemma \ref{key-lem} below) without using \v{C}ech representatives.\
\footnote{We were not able to understand Markarian's proof of a similar statement, \cite[Prop.\ 5]{Markarian}.} 

The proof of Theorem A uses analogs of some constructions of Markarian \cite{Markarian} for matrix factorizations.
Recall that he used the Atiyah classes to equip the shifted tangent bundle $T_X[-1]$ with a structure of
a Lie algebra in the derived category of $X$ (this construction goes back to Kapranov \cite{Kapranov}), 
and showed that in appropriate sense the universal enveloping of
this algebra can be identified with the sheafified Hochschild cohomology $\und{HH}^*(X)$. Furthermore,
the Hochschild-Kostant-Rosenberg isomorphism can be viewed as an analog of the PBW theorem in this case.

It turns out that there is a similar Lie context for the sheafified Hochschild cohomology of the category of matrix
factorizations, $\und{HH}^*\MF(X,w)$. The general principle is that the picture for $w\neq 0$ should be a deformation of the picture for $w=0$.
In Lie theory there is a well known way of deforming the universal enveloping algebras $U(\fg)$ of a Lie algebra starting from 
a central extension of Lie algebras
\begin{equation}\label{Lie-central-ext-eq}
0\to k\cdot\unit \to \wt{\fg}\to \fg\to 0.
\end{equation}
Namely, one can view $\unit\in \wt{\fg}$ as a central element of $U(\wt{\fg})$ and consider the quotient $U(\wt{\fg})/(\unit-1)$, which is a deformation of
$U(\wt{\fg})/(\unit)\simeq U(\fg)$.

Given a function $w$, we will equip the $\Z/2$-graded complex 
$$L_w:=[T_X\rTo{i_{dw}}\OO_X]$$ 
(where $\OO_X$ is placed in degree $0$) with a structure of a Lie algebra, so that
the exact triangle
$$\OO_X\to [T_X\rTo{i_{dw}} \OO_X]\to T_X[1]$$
can be viewed as a central extension of Lie algebras 
in $D(X)$, the $\Z/2$-graded derived category of $X$ (see Sec.\ \ref{general-Lie-sec}).
Note that such a construction wouldn't work in the usual $\Z$-graded derived category since it is $T_X[-1]$ that has the Lie algebra structure,
not $T_X[1]$.

Extending the picture of Markarian \cite{Markarian}, we will show that 
$\und{HH}^*\MF(X,w)$ can be viewed
as the corresponding quotient of the universal enveloping algebra, 
$U(L_w)/(\unit-1)$.

In Section \ref{graded-sec} we consider a version of the above picture for a $\Z$-graded category of matrix factorizations defined in 
the presence of a $\G_m$-action. More precisely, we assume that $X$ is equipped with a $\G_m$-action and we have a function $W$ on $X$
satisfying $W(\la x)=\la W(x)$. In this context one has a natural $\Z$-graded dg-category $\MF_{\G_m}(X,W)$ of {\it $\G_m$-equivariant matrix
factorizations of $W$}. We prove that the analog of Theorem A holds in this context. In the particular case when $W=0$ and the action of
$\G_m$ is trivial, the category $\MF_{\G_m}(X,0)$ is equivalent to the usual $\Z$-graded derived category of $X$, so in this
case we recover the classical picture for the latter category as described in  \cite{Caldararu-HKR}, \cite{Ram}.

Note that a more natural context for $\Z$-graded categories of matrix factorizations involves equivariant matrix factorizations with respect
to an algebraic group $\Ga$ equipped with a surjective homomorphism $\chi:\Ga\to \G_m$ such that $\ker(\chi)$ is finite. However,
the corresponding picture already has a stacky flavor. We intend to consider it elsewhere, along with a stacky version of
Theorem A.

\medskip

\noindent
{\it Conventions and notation}. We work with matrix factorizations of a regular function $w$ on a smooth scheme
$X$ of dimension $n$ over a field of characteristic $0$. By $\MF(w)=\MF(X,w)$ we denote the corresponding derived category of 
matrix factorizations.
Whenever we need to use the HKR-isomorphism for sheafified Hochschild
homology of the category of matrix factorizations of $w$, we assume that $0$ is the only critical value of $w$.
We denote by $D(X)$ the $\Z/2$-graded derived category of $X$, and by $D(X)_{\Z}$ the usual $\Z$-graded derived category.

\medskip

\noindent
{\it Acknowledgments}.
We thank the anonymous referee for useful comments.
B.K. is partially supported by the KIAS individual grant MG016403.
A.P. is partially supported by the NSF grant DMS-1700642 and within the framework of the HSE University Basic Research Program
by the Russian Academic Excellence Project `5-100'.
He also would like to thank the Korea Institute for Advanced Study and Institut des Hautes Etudes Scientifiques, where parts of this work were done, for hospitality
and excellent work conditions.

\section{Definition of the Atiyah class}

\subsection{Global definition}\label{global-atiyah-def-sec}
Recall that the Atiyah class for vector bundles (or bounded complexes thereof) is defined using
the canonical exact sequence
\begin{equation}\label{can-ext-X2-eq}
0\to \De_*\Om^1_X\to \OO_{\De^{(2)}}\to \OO_\De\to 0.
\end{equation}
Here $\OO_{\De^{(2)}}:=\OO_{X^2}/I_\De^2$, where
$I_\De\sub \OO_{X^2}$ is the ideal sheaf of the diagonal $\De\sub X^2$.
For a vector bundle $E$ (or any object in $D(X)$), we tensor the above sequence with $p_2^*E$ and apply $p_{1*}$ to
get the extension sequence 
$$0\to \Om^1_X\ot E\to J(E)\to E\to 0$$
representing $\at(E)$.

Equivalently, we can view $[ \De_*\Om^1_X  \to \OO_{\De^{(2)}}]$ (in degrees $-1$ and $0$) 
as a complex quasi-isomorphic to $\OO_\De$. Then the natural projection to $(\Om^1_X)_\De[1] := \De_*\Om^1_X [1]$ gives a morphism
in $D(X^2)$,
\begin{equation}\label{univ-atiyah-eq}
\at^{\univ}:\OO_\De\to(\Om^1_X)_\De[1],
\end{equation}
called the {\it universal Atiyah class}. From this morphism of kernels we get a morphism of functors,
whose value on $E$ is $\at(E)$.

\begin{lem}\label{univ-at-lem} 
The universal Atiyah class, $\at^{\univ}$, is equal to the composition
$$\OO_\De\rTo{\at(\OO_\De)}  \Om^1_{X^2}\ot \OO_\De[1]\simeq \De_*\De^*\Om^1_{X^2}[1]\to \De_*\Om^1_X[1],$$
where the last arrow is induced by the canonical map $\De^*\Om^1_{X^2}\to \Om^1_X$.
Equivalently, if we use the decomposition $\Om^1_{X^2}=p_1^*\Om^1_X\oplus p_2^*\Om^1_X$, then
$\at^{\univ}$ is equal to the component 
$$\at^1(\OO_\De)\in \Hom(\OO_\De,p_1^*\Om^1_X\ot \OO_\De[1])\simeq \Hom(\OO_\De,\De_*\Om^1_X[1])$$
of $\at(\OO_\De)$.
\end{lem}

\Pf . (See \cite[Sec.\ 5.5]{RW}.) The class $\at(\OO_\De)$ corresponds to the canonical extension
$$0\to \Om^1_{X^2}\ot\OO_{\De}\to J(\OO_\De)\to \OO_\De\to 0,$$
where $J(?)$ denotes the sheaf of 1st order jets. By definition, $J(\OO_\De)$ is obtained as the push-forward
of $\OO_{X^4}/(J_{\De_{13}}+J_{\De_{24}})^2\ot \OO_{\De_{34}}$ under $p_{12}:X^4\to X^2$, where $\De_{ij}\sub X^4$ 
are partial diagonals. Using the identification of $\De_{34}$ with $X^3$ we get 
$$J(\OO_\De)\simeq p_{12*}\OO_{X^3}/(J_{\De_{13}}+J_{\De_{23}})^2.$$ 
Hence, we have a natural map 
$$J(\OO_\De)\simeq p_{12*}\OO_{X^3}/(J_{\De_{13}}+J_{\De_{23}})^2\to p_{12*}\OO_{X^3}/(J_{\De_{13}}^2+J_{\De_{23}})\simeq \OO_{X^2}/J_{\De}^2.$$
One can check that it is compatible with the projection $\Om^1_{X^2}\ot\OO_{\De}\to p_1^*\Om^1_X\ot \OO_\De$.
Thus, we get a morphism of exact sequences
\begin{diagram}
0&\rTo{}& \Om^1_{X^2}\ot\OO_{\De}&\rTo{}&J(\OO_\De)&\rTo{}&\OO_\De&\rTo{}&0\\
&&\dTo{}&&\dTo{}&&\dTo{}\\
0&\rTo{}& p_1^*\Om^1_X\ot\OO_{\De}&\rTo{}&\OO_{\De^{(2)}}&\rTo{}&\OO_\De&\rTo{}&0
\end{diagram}
and our assertion follows.
\ed

We observe that the sequence \eqref{can-ext-X2-eq} has the following analog in the category of
(coherent) matrix factorizations $\MF(\wt{w})$, where $\wt{w}=w\ot 1-1\ot w\in H^0(X^2,\OO)$.
Note that we have a natural functor 
$$\De_*:\MF(X,0)\to \MF(X^2,\wt{w}),$$
since $\wt{w}|_\De=0$. 
We denote by $\OO_\De^{\wt{w}}\in  \MF(X^2,\wt{w})$ the image of $\OO_X$ under this functor.
Let us 
define the matrix factorizaton $\OO^{(2)}_{\De,\wt{w}}$ of $\wt{w}$ as follows:
$$\OO^{(2)}_{\De,\wt{w}}=\OO_{X^2}/I_\De^2\oplus \OO_\De[1]
%[\OO^{(2)}_{\De,\wt{w}}]^0=\OO_{\De^{(2)}}=\OO_{X^2}/I_\De^2, \ \ 
%[\OO^{(2)}_{\De,\wt{w}}]^1=\OO_\De,
$$
$$\de_0=-\wt{w}\mod I_{\De}^2=-dw:\OO_\De\to I_{\De}/I_{\De}^2\sub \OO_{\De^{(2)}}, \ \ 
\de_1=-1:\OO_{\De^{(2)}}\to\OO_\De.$$
Then we have an exact sequence of matrix factorizations
\begin{equation}\label{jet-kernel-mf-ex-seq}
0\to \De_*[\OO_X\rTo{dw} \Om^1_X][1]\to \OO^{(2)}_{\De,\wt{w}}\to \OO_\De^{\wt{w}}\to 0
\end{equation}
which can be viewed as a morphism in the derived category of $\MF(X^2, \wt{w})$,
$$\hat{\at}^{\univ}:\OO_\De^{\wt{w}}\to \De_*[\OO_X\rTo{dw} \Om^1_X].$$
Note that by the definition of the shift functor on complexes, the complex $\De_*[\OO_X\rTo{dw} \Om^1_X][1]$ has $\OO_X$ in odd degree,
$\Om^1_X$ in even degree and the differential is given by $-dw$. This is compatible with the sign in the definition of
the differential on $\OO^{(2)}_{\De,\wt{w}}$.

Now given a matrix factorization $E\in \MF(X, w)$, we can tensor the exact sequence \eqref{jet-kernel-mf-ex-seq} with $p_2^*E$ and then apply
$Rp_{1*}$. This will produce a class
$$\hat{\at}(E)\in \Hom^0(E,[\OO_X\rTo{dw} \Om^1_X]\ot E).$$

Equivalently, we can obtain $\hat{\at}(E)$ from the morphism of kernels
$$\OO_\De^{\wt{w}}\to \De_*[\OO_X\rTo{dw} \Om^1_X]$$
in $\MF(\wt{w})$ corresponding to the sequence \eqref{jet-kernel-mf-ex-seq}.
By definition, this morphism in $\MF(\wt{w})$ uses the quasi-isomorphism
of the cone of the first arrow in \eqref{jet-kernel-mf-ex-seq} with $\OO_{\De}^{\wt{w}}$.
It will be convenient to use a slightly more compact resolution of $\OO_{\De}^{\wt{w}}$,
which should be viewed as a curved analog of the resolution $[(\Om^1_X)_\De\to \OO_{\De^{(2)}}]$ of $\OO_\De$.
Namely, let us equip $\OO_{\De^{(2)}}\oplus (\Om^1_X)_\De[1]$ with the structure of a matrix factorization of $\wt{w}$
using the maps
$$\de_0:(\Om^1_X)_\De\simeq I_\De/I_\De^2\hra \OO_{\De^{(2)}}, \ \ 
\de_1:\OO_{\De^{(2)}}\rTo{1} \OO_{\De}\rTo{dw}(\Om^1_X)_\De .$$

\begin{lem}\label{mf-atiyah-basic-lem} 
(i) The natural map in $\MF(\wt{w})$,
$$q:[\OO_{\De^{(2)}}\oplus \De_*\Om^1_X[1],\de]\to \OO_\De^{\wt{w}},$$
induced by the projection $\OO_{\De^{(2)}}\to \OO_{\De}$,
is an isomorphism in the derived category 
and $\hat{\at}^{\univ}$ is equal to the composition of $q^{-1}$ with
the natural morphism in $\MF(\wt{w})$,
\begin{equation}\label{atiyah-res-map}
[\OO_{\De^{(2)}}\oplus \De_*\Om^1_X[1],\de]\to \De_*[\OO_X\rTo{dw}\Om^1_X],
\end{equation}
which is identity on $\De_*\Om_X^1[1]$ and the natural projection to $\OO_{\De}$ on $\OO_{\De^{(2)}}$.

\noindent
(ii) The composition 
$$\OO_{\De}^{\wt{w}}\to \De_*[\OO_X\rTo{dw}\Om_X^1]\to \De_*\OO_X=\OO_{\De}^{\wt{w}}$$
is the identity map. Hence,
the image of $\hat{\at}(E)$ under the natural projection
$\Hom^0(E,E\ot [\OO_X\rTo{dw}\Om^1_X])\to \Hom^0(E,E)$ is the identity element
$\id_E$.

\noindent
(iii) The universal Atiyah class $\hat{\at}^{\univ}$ is obtained as the composition of the Atiyah class 
$$\hat{\at}(\OO_{\De}^{\wt{w}}):\OO_\De^{\wt{w}}\to [\OO_{X^2} \rTo{d\wt{w}}\Om^1_{X^2}]\ot \OO_{\De}^{\wt{w}}$$
with the projection 
$$[\OO_{X^2}\rTo{d\wt{w}}\Om^1_{X^2}]\ot \OO_{\De}^{\wt{w}}\to
[\OO_{X^2} \rTo{d(w\ot 1)}p_1^*\Om^1_X]\ot \OO_{\De}^{\wt{w}}\simeq
\De_*[\OO_X\rTo{dw}\Om^1_X].$$
\end{lem}

\Pf . (i) Let $f:\De_*[\OO_X\rTo{dw} \Om^1_X][1]\to \OO^{(2)}_{\De,\wt{w}}$
be the map from the sequence \eqref{jet-kernel-mf-ex-seq} ($f$ is the identity on $\OO_{\De}$ and
is the natural embedding into $\OO_{\De^{(2)}}$ on $\De_*\Om^1_X$).
We have a natural quasi-isomorphism
$$\Cone(f)\to [\OO_{\De^{(2)}}\oplus (\Om^1_X)_\De[1],\de]$$
which is identical on $\De_*\Om_X^1[1]$ and $\OO_{\De^{(2)}}$,
zero on $\OO_\De[1]$ and equal to $-dw:\OO_\De\to \De_*\Om^1_X$ on $\OO_\De$.
One can check that its composition with \eqref{atiyah-res-map} is homotopic
to the canonical projection $\Cone(f)\to \De_*[\OO_X\rTo{dw} \Om^1_X]$
(one uses $\id:\OO_\De[1]\to \OO_\De[1]$ as a homotopy).
This implies our assertion.

\noindent
(ii) This follows easily from (i): the composition of \eqref{atiyah-res-map} with
the natural projection $\De_*[\OO_X\rTo{dw}\Om^1_X]\to \De_*\OO_X=\OO_{\De}^{\wt{w}}$
is exactly the map $q$.

\noindent
(iii) This is proved similarly to Lemma \ref{univ-at-lem}.
\ed

\subsection{\v{C}ech representative}

Let $(U_i)$ be an affine open covering of $X$. Over every $U_i$ we can choose an algebraic connection
$$\nabla_i:E|_{U_i}\to\Om^1_{U_i}\ot E_{U_i}$$
which is even, i.e., compatible with the $\Z/2$-grading on $E|_{U_i}$.
Over each intersection $U_{ij}=U_i\cap U_j$ we have a $1$-form with values in $\und{\End}^0(E)$, 
$\a_{ij}\in \Om^1\ot \und{\End}^0(E)|_{U_{ij}}$, such that 
$$\nabla_j-\nabla_i=\a_{ij}.$$

Assume for simplicity that $X$ is separated, so all intersections $U_{i_1\ldots i_k}$ are still affine. 
%(BK: The condition that $U_{i_1\ldots i_k}$ are still affine seems unnecessary.)
Then for any matrix factorizations of $w$, $E$ and $F$, we can calculate the space
$\Hom^0(E,F)$ as the $0$th cohomology of the $\Z/2$-graded complex
$$(C^\bullet(\und{\Hom}(E,F)),[\de,?]+d_C)$$
where $C^\bullet(?)$ denotes the \v{C}ech complex and $d_C$ is the \v{C}ech differential.
More precisely, the differential on $\a\in C^p(\und{\Hom}(E,F))$ is $(-1)^p[\de,\a]+d_C(\a)$.

\begin{prop}\label{cech-formula-prop} 
For a matrix factorization $E\in \MF(X, w)$, the Atiyah class
$\hat{\at}(E)$ is represented by the
cocycle
$$(\id_E,-[\nabla_i,\de],\a_{ij})\in C^\bullet([\OO\rTo{dw} \Om^1_X]\ot \und{\End}(E)).$$
\end{prop}

\Pf . We use the following general fact: if 
\begin{equation}\label{ex-seq-mf}
0\to E_1\rTo{f} E_2\rTo{g} E_3\to 0
\end{equation}
is an exact sequence of matrix factorizations then the corresponding class in $\Hom^1(E_3,E_1)$ is represented
by the following \v{C}ech cocycle. First, we find local retractions $r_i:E_2\to E_1$ of the embedding $E_1\to E_2$,
which are morphisms of $\Z/2$-graded $\OO$-modules over $U_i$.
Then we consider 
$\a_{ij}:E_3\to E_1$ over $U_{ij}$ such that
$\a_{ij}g=r_j|_{U_{ij}}-r_i|_{U_{ij}}$. 
On the other hand, we consider $\b_i:E_3\to E_1[1]$ over $U_i$ such that
$\b_ig=[\de,r_i]=\de_{E_1}r_i-r_i\de_{E_2}$. Now we claim that the \v{C}ech cocycle $c=(\b_i,\a_{ij})$ represents the 
class corresponding to our extension. Indeed, by definition, this class corresponds to the obvious
projection $\Cone(f)\to E_1[1]$ under the isomorphism 
$$\Hom^1(E_3,E_1)\rTo{\sim}\Hom^1(\Cone(f), E_1)$$
induced by $g$. Now the image of $c$ under the morphism  
$$C^1(\und{\Hom}(E_3,E_1))\to C^1(\und{\Hom}(\Cone(f),E_1))$$
is given by the cocycle 
$$(\b_ig,\a_{ij}g)=([\de,r_i],r_j-r_i).$$
Subtracting the coboundary of the element $(r_i)\in C^0(\und{\Hom}(\Cone(f),E_1))$, we get
the cocycle given by $(r_if)$, i.e., by $\id_{E_1}$, as claimed.

We apply the above general fact to the sequence
\begin{equation}\label{Jet-mf-seq}
0\to [\OO_X\rTo{dw}\Om^1_X]\ot E[1]\to \hat{J}(E)\to E\to 0
\end{equation}
obtained from \eqref{jet-kernel-mf-ex-seq} by tensoring with $p_2^*E$ and taking the push-forward $p_{1*}$.
Note that 
$$\hat{J}(E)=J(E)\oplus E[1]$$ 
as a $\Z/2$-graded vector bundle.
% and the differential on it has as components $J(\de)$, $\de_{E[1]}$,
%the natural projection $J(E)\to E$, and the embedding $E\rTo{dw} E\ot\Om^1_X\sub J(E)$.
A connection $\nabla_i$ on $E|_{U_i}$ can be viewed as a retraction 
$$\nabla_i:J(E)\to \Om^1_X\ot E: 1\ot s\mapsto \nabla_i(s).$$
This leads to the claimed formula.
%$$-\de_{\hat{J}(E)}\nabla_i+\nabla_i\de_E=(-J(\de_E)\nabla_i+\nabla_i\de(E),-\id_E)$$
\ed

\begin{rem}\label{trivial-bun-rem} 
In the case when the underlying vector bundle $E$ is trivialized, we can take $\nabla$ to be
the corresponding connection (defined globally). Then $[\nabla,\de]$ is the matrix of $1$-forms obtained by taking differentials of
the entries of $\de$ (viewed as a matrix of functions).
Using this one can check that in this case the formula of Theorem A for the Chern character is equivalent to the one
obtained in \cite{PV-Chern} in the case of matrix factorizations over a regular local ring.
\end{rem}

Using \v{C}ech representatives one can easily derive the following compatibility of the Atiyah class construction with
the tensor product of a matrix factorization by a complex of vector bundles.

\begin{lem}\label{atiyah-tensor-prod-lem}
For $E\in\MF(X,w)$ and $F\in D(X)$, let us consider $E\ot F\in\MF(X,w)$. Then one has
$$\hat{\at}(E\ot F)=\hat{\at}(E)\ot \id+\id\ot \at(F).$$ 
\end{lem}

\subsection{Dolbeault representative}

In the case when $X$ is a complex manifold and $w$ is a holomorphic function,
the space $\Hom^0(E,F)$ can be calculated using Dolbeault complex,
$$(\Om^{0,*}(\und{\Hom}(E,F)),[\de,?]+\dbar).$$
%More precisely, for $\a\ot f\in \Om^{0,p}(\und{\Hom}(E,F))$ the differential is given by $(-1)^p\a\ot [\de,f]+\dbar(

In the case $w=0$ it is well known that the Atiyah class is represented by the $(1,1)$-part of the curvature
of any $C^\infty$-connection on $E$, compatible with the holomorphic structure.

We have the following analog for matrix factorizations.

\begin{prop} Let $(E,\de)$ be a holomorphic matrix factorization of $w$. Let $\nabla$ be an even $C^\infty$-connection
on $E$, compatible with the holomorphic structure, and let $F^{1,1}$ be the $(1,1)$-part of the curvature of $\nabla$.
Then $\hat{\at}(E)$ is represented by the cocycle
$$(\id_E,-[\nabla,\de],F^{1,1})\in \Om^{\le 1,*}(\und{\End}(E,E)),$$
where the differential on the latter complex is given by $[\de,?]+\we dw+\dbar$.
\end{prop}

\Pf . The proof is similar to that of Proposition \ref{cech-formula-prop}, with the \v{C}ech resolution replaced 
by the Dolbeault resolution. First, one checks that for an exact sequence of matrix factorizations \eqref{ex-seq-mf},
a choice of $C^\infty$-retractions $r:E_2\to E_1$ gives a Dolbeault representative for the corresponding class in $\Hom^1(E_3,E_1)$.
Namely, one should consider 
$$(\b,\a)\in\Om^{0,0}(\und{\Hom}(E_3,E_1)_1)\oplus \Om^{(0,1)}(\und{\Hom}(E_3,E_1)_0),$$
where 
$$\a g=\dbar(r), \ \ \b g=[\de,r]$$
(the proof uses the Dolbeault complex of $\Cone(f)$, similarly to Proposition \ref{cech-formula-prop}).

Now we apply this for the exact sequence \eqref{Jet-mf-seq}.
We use $\nabla^{1,0}$, the $(1,0)$-part of the connection $\nabla$, to get a $C^\infty$-retraction of 
the embedding $\Om^{1,0}\ot E\to J(E)$. This leads to
the claimed formula, where the $(1,1)$-curvature $F^{1,1}$ appears as $[\dbar,\nabla^{1,0}]$.
\ed

\section{Lie algebra analogies}

\subsection{Lie bracket in the case $w=0$}
Recall that the Atiyah class $\at(\Om^1):\Om^1_X\to \Om^1_X\ot \Om^1_X[1]$ 
factors through a map $\Om^1_X\to S^2(\Om^1_X)[1]$ whose dual can be viewed as a Lie bracket on
$T_X[-1]$ (see \cite{Kapranov}, \cite{Markarian}).

In addition, there is a morphism 
\begin{equation}\label{iota-eq}
\iota:T_X[-1]\to \und{HH}^*(X)=Rp_{1*}\und{\End}(\OO_\De)
\end{equation}
obtained by adjunction from the component
$$\at^1(\OO_\De)\in \Hom(\OO_\De,p_1^*\Om_X^1\ot \OO_\De[1])\simeq\Hom(p_1^*T_X[-1],\und{\End}(\OO_\De))$$
of the Atiyah class $\at(\OO_\De)\in \Hom(\OO_\De,\Om_{X^2}^1\ot\OO_\De[1])$.

\subsection{HKR-isomorphisms in the case $w=0$}

Recall that in \cite[Theorem 1]{Markarian} Markarian showed that 
$\und{HH}^*(X)$ can be identifed with the universal enveloping algebra $U(T_X[-1])$. By definition, this
means that the above map $\iota$ satisfies the identity 
\begin{equation}\label{env-alg-identity}
\iota([x,y])=\iota(x)\iota(y)-\iota(y)\iota(x)
\end{equation}
understood as the equality of morphisms $T_X[-1]\ot T_X[-1]\to \und{HH}^*(X)$, and the morphism
\begin{equation}\label{I-coh-abs}
I^{abs}:S^*(T_X[-1])\to T^*(T_X[-1])\to \und{HH}^*(X),
\end{equation}
induced by $\iota$ and by the multiplication on $\und{HH}^*(X)$, is an isomorphism
(see \cite[Def.\ 4]{Markarian}).

One of the key tools in the arguments of \cite{Markarian} is the natural duality between
$\und{HH}_*(X)$ and $\und{HH}^*(X)$ induced by the canonical functional
$$\vareps:\und{HH}_*(X)=\De^*(\OO_\De)\to \OO_X,$$
and by the canonical action 
$$\DD_0:\und{HH}^*(X)\ot \und{HH}_*(X)\to  \und{HH}_*(X) .$$ 

\begin{lem}\label{pairing-0-lem} The composition
$$\vareps\circ\DD_0:\und{HH}^*(X)\ot \und{HH}_*(X)\to \OO_X$$
is a perfect pairing, which corresponds to the natural isomorphism
\begin{align*}
&\und{\Hom}(\und{HH}_*(X),\OO_X)=\und{\Hom}(\De^*\OO_\De,\OO_X)\simeq 
Rp_{1*}\De_*\und{\Hom}(\De^*\OO_{\De},\OO_X)\\
&\simeq Rp_{1*}\und{\Hom}(\OO_\De,\OO_\De) .
\end{align*}
\end{lem}

\Pf . Note that for any morphism $f:X\to Y$ and sheaves $F$ on $Y$ and $G$ on $X$, the composition
of the natural maps
$$\und{\Hom}(F,f_*G)\rTo{f^*} Rf_*\und{\Hom}(f^*F,f^*f_*G)\rTo{\can_G} Rf_*\und{\Hom}(f^*F,G)$$
is precisely the (sheafified) adjunction isomorphism.
Applying this to $f=\De$, $F=\OO_\De$, $G=\OO_X$, we obtain that the composition
$$\und{\Hom}(\OO_\De,\OO_\De)\rTo{\De^*}\De_*\und{\Hom}(\De^*\OO_\De,\De^*\OO_\De)\rTo{\eps}
\De_*\und{\Hom}(\De^*\OO_\De,\OO_X)$$
is the adjunction isomorphism. Now the assertion follows from the fact that applying $Rp_{1*}$ to the
first arrow we get the map $\und{HH}^*(X)\to \und{\Hom}(\De^*\OO_\De,\De^*\OO_\De)$ defining $\DD_0$.
%It is enough to prove this after passing to a completed local ring of a closed point in $X$, which we can
%identify with $k[[x_1,\ldots,x_n]]$. In this case we can use the diagonal Koszul matrix factorization to represent
%$\OO_\De^{\wt{w}}$???
\ed

Dualizing the composition
$$S^*(T_X[-1])\ot\und{HH}_*(X)\rTo{I^{abs}\ot\id}\und{HH}^*(X)\ot \und{HH}_*(X)\rTo{\vareps\circ \DD_0}\OO_X$$
we get a map
\begin{equation}\label{I-hom-abs}
I_{abs}:\und{HH}_*(X)\to S^*(T_X[-1])^\vee\simeq S^*(\Om_X[1]),
\end{equation}
which is an isomorphism (see \cite[Prop.\ 3]{Markarian}). 

Note that essentially by definition, our isomorphisms $I^{abs}$ and $I_{abs}$ are compatible with the duality of
Lemma \ref{pairing-0-lem} and the natural duality between $S^*(T_X[-1])$ and $S^*(\Om_X[1])$.

\begin{thm}\label{HKR-thm}
The isomorphisms $I^{abs}$ and $I_{abs}$ in $D(X)$ coincide with the Hochschild-Kostant-Rosenberg isomorphisms
$I^{HKR}$ and $I_{HKR}$,
given by the explicit chain maps in \cite{Caldararu-HKR}.
\end{thm}

%It is easy to see that the $i$th power of the universal Atiyah class by an explicit chain map from a resolution of $\OO_\De$
%given by
%???
%Then Caldararu constructs an explicit chain map from the representative complex of $\De^*\OO_\De$, obtained from the 
%completed bar-resolution $\BB_\bullet(X)$ of $\OO_\De$ to the above resolution of $\OO_\De$.

%\begin{rem} In the affine case Theorem \ref{HKR-thm} follows much easier ??? 
%\end{rem}

\Pf .
Recall that the action of $T_X[-1]$ on $\und{HH}_*(X)$ is obtained by applying
$\De^*$ to the map $p_1^*T_X[-1]\ot \OO_\De\to \OO_\De$ given by $\at^{\univ}=\at^1(\OO_\De)$.
We are going to realize the latter map by an explicit chain map of complexes on $X^2$,
replacing $\OO_\De$ by its completed bar-resolution $\BB_\bullet$ (see \cite{Yekut}, \cite[Sec.\ 4]{Caldararu-HKR}).
Recall that $\BB_i$ is the push-forward to $X^2$ of the formal completion of $X^{i+2}$ along the small diagonal.
We denote local sections of $\BB_i$ as $[a_0\ot a_1\ot\ldots\ot a_{i+1}]$, where $a_j$ are local functions on $X$.

We claim that the map $\at^{\univ}$ is represented by the chain map
$$p_1^*T_X\ot \BB_\bullet\to \BB_\bullet[1]: 
v\ot [a_0\ot a_1\ot\ldots\ot a_{i+1}]\mapsto (-1)^{i+1} [a_0 v(a_1)\ot a_2\ot\ldots \ot a_{i+1}].$$
%({\bf BK: Where does the sign $\pm$ come from?})
Indeed, it is enough to look at the induced map of complexes 
$p_1^*T_X\ot \BB_\bullet\to \OO_\De[1],$
or equivalently, $\BB_\bullet\to \De_*\Om^1_X[1]$
induced by the map $f_1:\BB_1\to \De_*\Om^1: [a_0\ot a_1\ot a_2]\mapsto a_2a_0da_1$.
Now we observe that there is a quasi-isomorphism of complexes,
\begin{diagram}
\BB_2&\rTo{}&\BB_1&\rTo{}&\BB_0\\
\dTo{}&&\dTo{f_1}&&\dTo{1}\\
0&\rTo{}&\De_*\Om^1_X&\rTo{}&\OO_{\De^{(2)}}
\end{diagram}
Since $\at^{\univ}$ is induced by the exact sequence \eqref{can-ext-X2-eq}, this implies our claim.
% that the chain map $\wt{\nu}$ realizes the action of $T_X[-1]$ on $\und{HH}_*(X)$. 

It follows that the map
$$(\at^{\univ})^{(i)}:p_1^*{\bigwedge}^i T_X\ot \OO_\De\to \OO_\De[i]$$
obtained by iterating $\at^{\univ}$ and restricting to skew-symmetric tensors, is represented by the chain map
$$p_1^*{\bigwedge}^i T_X\ot \BB_\bullet\to \BB_\bullet[i]: 
(v_1\we\ldots\we v_i)\ot [a_0\ot a_1\ot\ldots\ot a_m]\mapsto (-1)^{im}[a_0 \lan v_1\we\ldots\we v_i, da_1\we\ldots\we da_i\ran\ot 
a_{i+1}\ot\ldots \ot a_m].$$
Composing with the projection $\BB_\bullet[i]\to \OO_\De[i]$, we see that
$(\at^{\univ})^{(i)}$ is represented by the chain map
$p_1^*{\bigwedge}^i T_X\ot \BB_\bullet\to \OO_\De[i]$ 
which corresponds by duality to the map
$$\BB^i\to p_1^*\Om^i_X\ot \OO_\De\simeq \De_*\Om^i_X: [a_0\ot a_1\ot\ldots\ot a_{i+1}]\mapsto 
a_{i+1}a_0da_1\we\ldots\we da_i$$
which is exactly the $i$th component of $I_{HKR}$. This proves the equality $I_{abs}=I_{HKR}$.
%On the other hand, one can easily check that the adjunction map $\eps:\und{HH}_*(X)=\De^*\OO_\De\to \OO_X$, 
%is represented by the chain map of complexes $\De^*\BB^\bullet\to \OO_X$, corresponding to the
%natural map $\De^*\BB^0\to \OO_X$.

%By Lemma \ref{I-abs-at-connection-lem}, 
%this amounts to checking that
%the composition
%$$\und{HH}_*(X)\rTo{I_{HKR}} S^*(\Om_X^1[1])\to \Om_X^1[1]$$
%is equal to the map $\De^*\OO_\De\to \Om_X^1[1]$ coming from
%the exact sequence \eqref{can-ext-X2-eq}.
%Thus, it is enough 
%to construct a chain map of complexes on $X^2$,

%such that $f_1$ agrees with the component of $I_{HKR}$.
%The map $f_1$ is given by $a_0\ot a_1\ot a_2\mapsto a_0a_2da_1$ and the required properties are easily checked.

Recall (see \cite[Sec.\ 4]{Caldararu-HKR}) that the other HKR-map
$$I^{HKR}:\bigoplus {\bigwedge}^iT_X[-i]\to Rp_{1*}\und{\End}(\OO_\De)$$
is essentially the dual of $I_{HKR}$: it is given by the composition
$$\bigoplus {\bigwedge}^iT_X[-i]=\und{\Hom}(\bigoplus\Om^i_X[i],\OO_X)\to \und{\Hom}(\De^*\OO_\De,\OO_X)\simeq
Rp_{1*}(\und{\Hom}(\OO_\De,\OO_\De),$$
where the last isomorphism follows from the adjunction of $(\De^*,\De_*)$.

On the other hand, it is easy to see that the composition
$$\und{HH}^*\ot \und{HH}_*\rTo{\DD_0} \und{HH}_*\rTo{\vareps}\OO_X$$
corresponds to the same duality, so the equality $I^{abs}=I^{HKR}$ follows from the equality $I_{abs}=I_{HKR}$.
\ed

\begin{rem} The above theorem can be easily deduced from the arguments in the proof of \cite[Prop.\ 4.4]{Caldararu-HKR}.
Note that $I^{HKR}$ (and hence $I^{abs}$) is not an algebra homomorphism with respect to the
natural algebra structures on $S^*(T_X[-1])$ and $\und{HH}^*(X)$; to become one it has to be twisted by the square root
of the Todd class (see \cite{Caldararu-HKR}, \cite{CVdB}).
\end{rem} 

%One application of Lemma \ref{derivation-0-lem} is the commutativity of the following diagram
%\begin{equation}\label{deriv-exp}
%\begin{diagram}
%T_X[-1]\ot E\ot\und{HH}_*(X)&\rTo{\id\ot\exp(\at(E))}& T_X[-1]\ot E\ot\und{HH}_*(X)\\
%\dTo{\DD_0+\at(E)\ot\id}&&\dTo{\DD_0}\\
%E\ot\und{HH}_*(X)&\rTo{\exp(\at(E))}& E\ot\und{HH}_*(X)
%\end{diagram}
%\end{equation}
%Is this true???
%which expresses the formula $D(\e

% and the observation that the corresponding
%action of the Lie algebra $T_X[-1]$ on $\und{HH}_*(X)$, as on the dual to its enveloping algebra,
%is induced (by applying $\De^*$) by the component of the Atiyah class, $\at^1(\OO_\De)$, viewed as a map
%$$p_1^*T_X[-1]\to\und{\End}(\OO_\De).$$

\subsection{The general case}\label{general-Lie-sec}

The (usual) Atiyah class of the complex $L_w^\vee=[\OO_X\rTo{dw}\Om^1_X]$ is an element
$$\at(L_w^\vee): L_w^\vee\to \Om^1_X[1]\ot L_w^\vee.$$
In the $\Z/2$-graded derived category we have a natural morphism $\Om^1_X[1]\to L_w^\vee$.
Thus, composing the map $\at(L_w^\vee)$ with this morphism, we get a map
$$L_w^\vee\to L_w^\vee\ot L_w^\vee,$$
or dualizing, a map
$$[\cdot,\cdot]: L_w\ot L_w\to L_w,$$
which factors through $T_X[1]\ot L_w$.

%The following lemma will show that our bracket is skew-symmetric and that
%$\OO_X\sub L_w$ is central with respect to it.

\begin{lem} The dual of the bracket $[\cdot,\cdot]$ factors in the $\Z/2$-graded derived category as a composition
$$L_w^\vee\to S^2\Om^1_X\to \Om^1_X\ot \Om^1_X\to L_w^\vee\ot L_w^{\vee}.$$
Hence, the bracket $[\cdot,\cdot]$ is skew-symmetric and $\OO_X\sub L_w$ is central with respect to it.
\end{lem}

\Pf . The proof is a slight variation of the proof of \cite[Prop.\ 1.1]{Markarian}.
By definition, the map $\at(L_w^\vee)$ corresponds to an exact sequence
of complexes
$$0\to p_{1*}[I_\De/I_\De^2\ot p_2^*L_w^\vee]\to p_{1*}[\OO_{X^2}/I_\De^2\ot p_2^*L_w^\vee]\to p_{1*}[\OO_\De\ot p_2^*L_w^\vee]\to 0.$$
Now we use the natural morphism induced by the canonical differential $d_{X^2}:\OO_{X^2}\to \Om^1_{X^2}$:
$$I_\De/I_\De^3\to \OO_{X^2}/I_\De^3\rTo{-d_{X^2}} \OO_{X^2}/I_\De^2\ot \Om^1_{X^2}\to \Om_{X^2}/I_\De^2\ot p_2^*\Om^1_X.$$
Note that under this morphism $f\ot 1-1\ot f\mod I_\De^3$ is sent to $1\ot p_2^*(df)$.
This morphism extends to a chain map of complexes
$$[\OO_X\rTo{\de} p_{1*}(I_\De/I_\De^3)]\to p_{1*}[\OO_{X^2}/I_\De^2\ot p_2^*L_w^\vee],$$
where $\de(f)=f\cdot (1\ot w-w\ot 1)$, and the map $\OO_X\to p_{1*}(\OO_{X^2}/I_\De^2)$ is given by $f\mapsto f\ot 1\mod I_\De^2$. 
Furthemore, this chain map extends naturally to a morphism of exact sequences of complexes
\begin{diagram}
0&\rTo{}&p_{1*}(I_\De^2/I_\De^3)[-1]&\rTo{}& [\OO_X\to p_{1*}(I_\De/I_\De^3)]&\rTo{}& [\OO_X\to p_{1*}(I_\De/I_\De^2)]&\rTo{}&0\\
&&\dTo{}&&\dTo{}&&\dTo{\sim}\\
0&\rTo{}&p_{1*}[I_\De/I_\De^2\ot p_2^*L_w^\vee]&\rTo{}&p_{1*}[\OO_{X^2}/I_\De^2\ot p_2^*L_w^\vee]&\rTo{}&p_{1*}[\OO_\De\ot p_2^*L_w^\vee]&\rTo{}&0
\end{diagram}
in which the leftmost vertical arrow can be identified with the natural map $S^2\Om^1_X[-1]\to \Om^1_X\ot L_w^\vee$. 
This implies our assertion.
\ed

By analogy with morphism \eqref{iota-eq} we want to define a morphism
\begin{equation}\label{iota-w-eq}
\iota_w:L_w\to \und{HH}^*(\MF(X,w))=Rp_{1*}\und{\End}(\OO_\De^{\wt{w}})
\end{equation}
in $D(X)$. 
For this we consider the universal Atiyah class
%$$\hat{\at}(\OO_\De^{\wt{w}})\in 
%\Hom(\OO_\De^{\wt{w}},\OO_\De^{\wt{w}}\ot (\OO_{X^2} \rTo{d\wt{w}} \Om_{X^2}^1)),$$
%and use the natural morphism
%$$[\OO_{X^2}\rTo{d\wt{w}} \Om_{X^2}^1)]\to p_1^*[\OO_X\rTo{dw}\Om_X^1]=p_1^*L_w^\vee$$
%to get an element
\begin{equation}\label{at1-diag-eq}
\hat{\at}^{\univ}=\hat{\at}^1(\OO_\De^{\wt{w}})\in \Hom(\OO_\De^{\wt{w}},p_1^*L_w^\vee\ot \OO_\De^{\wt{w}})
\end{equation}
(see Lemma \ref{mf-atiyah-basic-lem}(iii)). Dualizing it we get a morphism
$$p_1^*L_w\to \und{\End}(\OO_\De^{\wt{w}})$$
from which $\iota_w$ is obtained by adjunction.

Note that to identify an associative algebra $U$ with the algebra of the form $U(\wt{\fg})/(\unit-1)$ for a central extension
of Lie algebras \eqref{Lie-central-ext-eq}, one has to provide a linear map $\iota:\wt{\fg}\to U$ satisfying 
the universal enveloping algebra identity \eqref{env-alg-identity}, such that $\iota(\unit)=1$ and 
the natural map 
$$S(\wt{\fg})/(\unit-1)\to U(\wt{\fg})/(\unit-1),$$
induced by $\iota$ is an isomorphism. Note the source of this map can be identified with $\varinjlim_i S^i(\wt{\fg})$,
which is a better expression for us since it makes sense also in non-abelian categories.

We want to check analogs of these properties for the map \eqref{iota-w-eq} in the derived category.
Similarly to the proof of \cite[Thm.\ 1]{Markarian}, the universal enveloping algebra identity
\eqref{env-alg-identity} for $\iota_w$ is equivalent to the condition that the skew-symmetrization of the composition
$$\OO_\De^{\wt{w}}\rTo{\hat{\at}^1(\OO_\De^{\wt{w}})}p_1^*L_w^\vee\ot \OO_\De^{\wt{w}}\rTo{\id\ot \hat{\at}^1(\OO_\De^{\wt{w}})}
p_1^*L_w^\vee\ot p_1^*L_w^\vee\ot \OO_\De^{\wt{w}}$$
is equal to the composition
$$\OO_\De^{\wt{w}}\rTo{\hat{\at}^1(\OO_\De^{\wt{w}})}p_1^*L_w^\vee\ot \OO_\De^{\wt{w}}\rTo{\at^1(p_1^*L_w^\vee)\ot \id}
p_1^*L_w^\vee\ot p_1^*L_w^\vee\ot \OO_\De^{\wt{w}}.$$
But this follows immediately from the commutative
diagram
\begin{diagram}
\OO_\De^{\wt{w}} &\rTo{\hat{\at}^1(\OO_\De^{\wt{w}})}&p_1^*L_w^\vee\ot \OO_\De^{\wt{w}}\\
\dTo{\hat{\at}^1(\OO_\De^{\wt{w}})} && \dTo{\hat{\at}^1 (p_1^*L_w^\vee\ot \OO_\De^{\wt{w}})}\\
p_1^*L_w^\vee\ot \OO_\De^{\wt{w}}&\rTo{\id\ot\hat{\at}^1(\OO_\De^{\wt{w}})}&
p_1^*L_w^\vee\ot p_1^*L_w^\vee\ot \OO_\De^{\wt{w}}
\end{diagram}
by applying Lemma \ref{atiyah-tensor-prod-lem} to the expression $\hat{\at}^1 (p_1^*L_w^\vee\ot \OO_\De^{\wt{w}})$.
%\si_{23}\circ\hat{\at}^1(\OO_\De^{\wt{w}})\ot\id
%\footnote{BK: Isn't $\hat{\at} (\OO_\De^{\wt{w}}\ot p_1^*L_w^\vee)$ instead of $\si_{23}\circ\hat{\at}^1(\OO_\De^{\wt{w}})\ot\id$ ?  }

The fact that $\hat{\at}(E)$ projects to the identity $\id_E$ implies (by taking $E=\OO_\De^{\wt{w}}$)
 that the composition 
$$\OO_X \rTo{\unit} L_w\rTo{\iota_w} \und{HH}^*(\MF(X,w))$$
is the natural embeding of a unit.

Lastly, we need to check that the map
\begin{equation}\label{I-abs-coh-w-eq}
I^{abs,w}:\varinjlim_i S^i(L_w)\to \und{HH}^*(\MF(X,w)),
\end{equation}
induced by $\iota_w$, is an isomorphism. 
It is easy to see that the limit here stabilizes and we have
$$\varinjlim_i S^i(L_w)=S^n(L_w)\simeq ({\bigwedge}^\bullet(T_X),i_{dw}),$$
where $n=\dim X$. Here the second isomorphism is dual to \eqref{Sn-isom}.
The fact that the map \eqref{I-abs-coh-w-eq} is an isomorphism can be checked formally locally
using the Koszul resolution of the diagonal matrix factorization (it also follows from Theorem \ref{HKR-w-thm} below
and from the results of \cite{LP}).

Similarly to the case of Lie algebras, where the PBW-theorem can be used to derive the Jacobi identity (see e.g.,
\cite[Ch.\ 5]{PP}), one can show that the properties proved above imply that the bracket $[\cdot,\cdot]$ on $L_w$
satisfies the Jacobi identity (we will not use this fact). 
%({\bf BK}: Why?)

Similarly to the case $w=0$, we have a canonical functional
$$\vareps:\und{HH}_*(\MF(X,w))=\De^*(\OO_\De^{\wt{w}})\to \OO_X,$$
coming from the adjoint pair of functors $(\De^*,\De_*)$ between $\MF(X,0)$ and $\MF(X^2,\wt{w})$.
On the other hand, we have the natural action of $\und{HH}^*(\MF(X,w))=Rp_{1*}\und{\End}(\OO_\De^{\wt{w}})$ on 
$\und{HH}_*(\MF(X,w))$,
$$\DD: \und{HH}^*(\MF(X,w))\ot \und{HH}_*(\MF(X,w))\to \und{HH}_*(\MF(X,w)).$$

\begin{lem}\label{pairing-lem} The composition
$$\vareps\circ\DD:\und{HH}^*(\MF(X,w))\ot \und{HH}_*(\MF(X,w))\to \OO_X$$
is a perfect pairing, which corresponds to the natural isomorphism
\begin{align*}
&\und{\Hom}(\und{HH}_*(\MF(X,w)),\OO_X)=\und{\Hom}(\De^*\OO_\De^{\wt{w}},\OO_X)\simeq 
Rp_{1*}\De_*\und{\Hom}(\De^*\OO_{\De}^{\wt{w}},\OO_X)\\
&\simeq Rp_{1*}\und{\Hom}(\OO_\De^{\wt{w}},\OO_\De^{\wt{w}}),
\end{align*}
induced by the adjoint pair of functors $(\De^*,\De_*)$ between $\MF(X,0)$ and $\MF(X^2,\wt{w})$.
\end{lem}

\Pf . The proof is very similar to that of Lemma \ref{pairing-0-lem} and is left to the reader. 
%It is enough to prove this after passing to a completed local ring of a closed point in $X$, which we can
%identify with $k[[x_1,\ldots,x_n]]$. In this case we can use the diagonal Koszul matrix factorization to represent
%$\OO_\De^{\wt{w}}$???
\ed
%and the natural action of $L_w$ on $\und{HH}_*(\MF(X,w))=\De^*(\OO_\De^{\wt{w}})$ is obtained by applying $\De^*$ to
%the Atiyah class component \eqref{at1-diag-eq}, viewed as a morphism
%$$p_1^*L_w\to\und{\End}(\OO_\De^{\wt{w}}).$$

Let us consider the composition
$$S^n(L_w)\ot \und{HH}_*(\MF(X,w))\rTo{I^{abs,w}\ot\id}\und{HH}^*(\MF(X,w))\ot \und{HH}_*(\MF(X,w))
\rTo{\vareps\circ\DD} \OO_X.$$
Dually, we get a morphism
\begin{equation}\label{I-abs-hom-w-eq}
I_{abs,w}:\und{HH}_*(\MF(X,w))\to S^n(L_w)^\vee\simeq [\Om^\bullet_X,\we dw],
\end{equation}
where the last isomorphism is \eqref{Sn-isom}.

\begin{thm}\label{HKR-w-thm} 
The maps $I^{abs,w}$ and $I_{abs,w}$ in $D(X)$ coincide with the maps $I_{HKR,w}$ and $I^{HKR,w}$
defined using the completed bar-resolution in \cite{LP} and \cite{Platt}.
\end{thm}

\Pf .
As in the proof of Theorem \ref{HKR-thm}, to check the equality $I_{abs,w}=I_{HKR,w}$,
we first realize 
$$\hat{\at}^{\univ}:p_1^*L_w\ot \OO_\De^{\wt{w}}\to \OO_\De^{\wt{w}}$$ 
by a closed map of
matrix factorizations.

Recall that the completed bar-resolution $(\BB_\bullet,b)$ is equipped with the second differential
$$B_w[a_0\ot a_1\ot\ldots\ot a_{m+1}]=\sum_{i=0}^m (-1)^i a_0\ot\ldots\ot a_i\ot w\ot a_{i+1}\ot \ldots\ot a_{m+1},$$
so that $\BB^{\wt{w}}:=(\BB_\bullet,b+B_w)$ is a (quasicoherent) matrix factorization of $\wt{w}=w\ot 1-1\ot w$.

Similarly to the case $w=0$, we define a closed morphism of matrix factorizations
$$p_1^*L_w\ot \BB^{\wt{w}}\to \BB^{\wt{w}}: (v+f)\ot [a_0\ot a_1\ot\ldots \ot a_m]=
(-1)^m[a_0v(a_1)\ot a_2\ot\ldots\ot a_m]+ [fa_0\ot a_1\ot\ldots\ot a_m].$$
We claim that it represents $\hat{\at}^{\univ}$.
Indeed, it is enough to consider the composed map
$$p_1^*L_w\ot \BB^{\wt{w}}\to \OO_\De^{\wt{w}}$$
and compare its dualization
\begin{equation}\label{at-w-resol}
\BB^{\wt{w}}\to p_1^*L_w^\vee\ot \OO_\De^{\wt{w}}\simeq \De_*[\OO_X\rTo{dw}\Om_X^1]
\end{equation} 
with $\hat{\at}^{\univ}$. It remains to observe that \eqref{at-w-resol}
factors as the composition
$$\BB^{\wt{w}}\to (\OO_{\De^{(2)}}\oplus(\Om_X^1)_\De[1],\de)\to (\OO_\De\oplus (\Om_X^1)_{\De}[1],dw)$$
where the first map is an isomorphism of resolutions of $\OO^{\wt{w}}_\De$, while the second
map induces $\hat{\at}^{\univ}$ by Lemma \ref{mf-atiyah-basic-lem}(i).

Considering the induced map $p_1^*S^nL_w\ot \BB^{\wt{w}}\to \BB^{\wt{w}}\to \OO_\De^{\wt{w}}$
we deduce that $I_{abs}$ coincides with
the Hochschild-Kostant-Rosenberg isomorphism $I_{HRK,w}$ given by the map
$$\De^*\BB^{\wt{w}}\to (\Om^\bullet_X,\we dw):[a_0\ot\ldots \ot a_{i+1}]\mapsto a_{i+1}a_0da_1\we\ldots\we da_i.$$

%This proves the equality $I_{abs,w}=I_{HKR,w}$. 
The equality $I^{abs,w}=I^{HKR,w}$ follows by duality as in the case $w=0$.
\ed

\section{Boundary-bulk map}

\subsection{Generalities}

Recall that the diagonal matrix factorization $\OO_\De^{\wt{w}}\in\MF(X^2,\wt{w})$
corresponds to the identity functor on $\MF(X,w)$.
The categorical trace functor can be identified with the composition
$$\Tr: \MF(X^2,\wt{w})\rTo{\De^*}D(X)\rTo{R\Ga}D(k).$$

Thus, the Hochschild homology of the category $\MF(X,w)$ can be computed as
$$HH_*(\MF(X,w))=R\Ga(X,\und{HH}_*(\MF(X,w)),$$
where
$$\und{HH}_*(\MF(X,w))=\De^*\OO_\De^{\wt{w}}.$$

Furthermore, we have an isomorphism
$$\und{HH}_*(\MF(X,w))\simeq [\Om^\bullet,\we dw].$$

The sheafified boundary-bulk map 
\begin{equation}\label{sh-bb-map}
\und{\End}(E)\to \und{HH}_*(\MF(X, w)),
\end{equation}
which is a map in $D(X)$,
is obtained by applying $\De^*$ to the evaluation morphism in $\MF(X^2,\wt{w})$,
$$\ev_E: E\boxtimes E^\vee \to \OO_\De^{\wt{w}}.$$
The latter morphism is obtained by dualization from the morphism
$$\eta_E:p_1^*E\to E_\De\simeq p_2^*E\ot  \OO_\De^{\wt{w}}$$ 
in $\MF(X^2,w\ot 1)$, which corresponds by adjunction to the isomorphism
$E\to Rp_{1*}(p_2^*E\ot \OO_\De^{\wt{w}})$.

Since the boundary-bulk map is obtained from $\ev_E$ by applying the categorical trace functor $\Tr$,
it is obtained from the sheafified boundary-bulk map \eqref{sh-bb-map} by passing to derived global sections.

\subsection{Exponentials}

Our exponentials are analogs of the following Lie theoretic construction. 
Let $\fg$ be a Lie algebra, $M$ a $\fg$-module.
Then for every $i\ge 0$, we have a morphism given by the iterated action,
$$u_i:\fg^{\ot i}\ot M\to M: x_1\ot\ldots\ot x_i\ot m\mapsto x_1\cdot(\ldots \cdot (x_{i-1} \cdot(x_i\cdot m))\ldots).$$
We denote by $s_i:S^i(\fg)\ot M\to M$ the restriction of $u_i$ to symmetric tensors.
We can think of $s_i$ as an element of $S^i(\fg)^*\ot\End(M)$.
Combining these elements together we get the element
$$\exp_M\in S^\bullet(\fg)^*\ot\End(M).$$

%Let us set $\AA=\bigoplus_{i\ge 0}\Om_X^i[i]$. This is a (super)commutative
%algebra object in $D(X)$. This induces an algebra structure on
%$\Hom(\OO_\De,\AA_\De)$. Namely, for $f,g\in\Hom(\OO_\De,\AA_\De)$ we define the product $fg$
%as the composition 
%$$\OO_\De\rTo{g} \De_*\AA\simeq p_1^*\AA\ot \OO_\De\rTo{\id\ot f}p_1^*\AA\ot \AA_\De\simeq
%(\AA\ot\AA)_\De\rTo{\mu}\AA_\De,$$
%where $\mu$ denotes the product on $\AA$. It is easy to check that it is associative.

%\begin{lem}\label{comm-lem}
%The subalgebra of $\Hom(\OO_\De,\AA_\De)$ generated by $\Hom(\OO_\De,\Om_X^1[1])$ over
%$\Hom(\OO_\De,\OO_\De)$ is commutative.
%\end{lem}

Now given an object $E$ in $D(X)$, the Atiyah class of $E$ defines a map
$$T_X[-1]\ot E\to E$$
which is an action of $T_X[-1]$ (viewed as a Lie algebra) on $E$.
Thus, we get the corresponding element
$$\exp_E\in \Hom(E,S^\bullet(\Om_X[1])\ot E)=\Hom(E,\bigoplus_i \Om^i_X[i]\ot E).$$
Unraveling the definitions, we see that
$$\exp_E=\exp(\at(E)),$$
where the right-hand side is defined in the standard way (see e.g., \cite[Sec.\ 4]{Caldararu-HKR}).

Similarly, we  can consider $\OO_\De$ as a module over $p_1^*T_X[-1]$ using the universal
Atiyah class, $\at^{\univ}:p_1^*T_X[-1]\ot \OO_\De\to \OO_\De$. This gives rise to the element
\begin{equation}\label{exp-at-univ-w0-eq}
\exp(\at^{\univ})\in \Hom(\OO_\De,\bigoplus_i p_1^*\Om^i_X[i]\ot \OO_\De)=\Hom(\OO_\De,\De_*\bigoplus_i \Om^i_X[i]).
\end{equation}
%We define $\exp(\at^{\univ})\in \Hom(\OO_\De,\AA_\De)$, applying the above algebra structure on
%$\Hom(\OO_\De,\AA_\De)$. We can formulate an analog of the differential equation for the exponential
%using the natural action of $T_X[1]$ on $\AA$.

It is easy to check that if we view $\exp(\at^{\univ})$ as a morphism of Fourier-Mukai kernels then the induced morphism
of functors $D(X)\to D(X)$, 
$$E\mapsto \bigoplus_i \Om^i_X[i]\ot E,$$
is precisely $\exp(\at(E))$.

In the same way for a matrix factorization $E\in\MF(X, w)$ we define 
$$\exp(\at(E))\in \Hom^0(E,(\Om_X^\bullet,\we dw)\ot E)$$
using the action of $L_w$ on $E$ given by $\hat{\at}(E)$, passing to the induced map
$$S^n(L_w)\ot E\to E$$
(where $n=\dim X$), dualizing, and using the isomorphism of $S^n(L_w^\vee)$ with $(\Om_X^\bullet,\we dw)$
(see \eqref{Sn-isom}). This is equivalent to the definition given in Theorem A.

Similarly, the universal Atiyah class,
$$\hat{\at}^{\univ}\in \Hom(\OO_{\De}^{\wt{w}},p_1^*L_w^\vee\ot\OO_{\De}^{\wt{w}})$$
gives an action of $p_1^*L_w$ on $\OO_\De^{\wt{w}}$
%which should be thought of as an analog of $1+\at^{\univ}$.
which induces an element
$$\exp(\at^{\univ})\in \Hom(\OO_{\De}^{\wt{w}},p_1^*S^n(L_w^\vee)\ot\OO_{\De}^{\wt{w}})\simeq 
p_1^*(\Om_X^\bullet,\we dw)\ot \OO_\De^{\wt{w}}.$$  

Again, if we view this as a morphism of Fourier-Mukai kernels then the induced morphism of functors $\MF(X,w)\to \MF(X,w)$
is given by $\exp(\at(E))$ on $E\in \MF(X,w)$.
%Let us set $\AA_w=(\Om_X^\bullet,\we dw)$. 

\subsection{Key lemma}

First, let us formulate the key assertion about the exponential of the universal Atiyah class \eqref{exp-at-univ-w0-eq} in the case $w=0$.
%Recall that we have the morphism
%$$\exp(\at^{\univ}):\OO_\De\to \De_*\bigoplus_i \Om_X^i[i].$$
%We have the following key fact.

\begin{lem}\label{key-lem}
One has a commutative triangle
\begin{diagram}
\OO_\De &\rTo{\can}& \De_*\De^*\OO_\De\\
&\rdTo{\exp(\at^{\univ})}&\dTo_{\De_*I_{abs}}\\
&& \De_*\bigoplus_i\Om_X^i[i]
\end{diagram}
\end{lem}

\Pf . Since the action of $T_X[-1]$ on $\und{HH}_*(X)=\De^*\OO_\De$ is obtained by applying $\De^*$ to
$\at^{\univ}$, by naturality of $\can$, we get the following commutative diagram
\begin{diagram}
p_1^*T_X[-1]\ot \OO_\De&\rTo{\id\ot\can}&p_1^*T_X[-1]\ot \De_*\De^*\OO_\De\\
\dTo{\at^{\univ}}&&\dTo{}\\
\OO_\De&\rTo{\can}&\De_*\De^*\OO_\De
\end{diagram}
where the right vertical arrow corresponds to the action of $T_X[-1]$ on $\und{HH}_*(X)$.
In other words, the map $\can:\OO_\De\to \De_*\De^*\OO_\De$ is compatible with the action of $p_1^*T_X[-1]$.
Hence, it is also compatible with the iteration of this action and its restriction to (skew)-symmetric tensors:
\begin{diagram}
p_1^*S^\bullet(T_X[-1])\ot \OO_\De&\rTo{\id\ot\can}&p_1^*S^\bullet(T_X[-1])\ot \De_*\De^*\OO_\De\\
\dTo{\exp(\at^{\univ})}&&\dTo{}\\
\OO_\De&\rTo{\can}&\De_*\De^*\OO_\De
\end{diagram}
where the left vertical arrow corresponds to $\exp(\at^{\univ})$ by duality.
Composing with the map $\De_*\vareps:\De_*\De^*\OO_\De\to \OO_\De$ (whose composition with $\can$ is
the identity), 
we get the commutative triangle
\begin{diagram}
p_1^*S^\bullet(T_X[-1])\ot \OO_\De&\rTo{\id\ot\can}&p_1^*S^\bullet(T_X[-1])\ot \De_*\De^*\OO_\De\\
&\rdTo{\exp(\at^{\univ})}&\dTo{}\\
&&\OO_\De
\end{diagram}
where the vertical arrow is the push-forward by $\De$
of the composition 
$$\vareps\circ \DD_0\circ (I^{abs}\ot \id):S^\bullet(T_X[-1])\ot \De^*\OO_\De\to \OO_X.$$ 
Now the assertion follows from the fact that $I_{abs}$ is obtained from
the latter map by dualization.
\ed

\begin{rem} Note that modulo Theorem \ref{HKR-thm} the assertion of Lemma \ref{key-lem}
is equivalent to that of Proposition 4.4 in \cite{Caldararu-HKR} (which refers to the standard HKR-isomorphism
defined using the completed bar-resolution) . Thus, we get a more conceptual proof
of that Proposition.
%with $I_{abs}$ replaced by , is exactly . However, the only place where
%the computations using global resolutions are needed is Theorem \ref{HKR-thm}.
\end{rem}

Now let us consider the case of matrix factorizations.

\begin{lem}\label{key-w-lem}
One has a commutative triangle
\begin{diagram}
\OO_\De^{\wt{w}} &\rTo{\can}& \De_*\De^*\OO_\De^{\wt{w}}\\
&\rdTo{\exp(\at^{\univ})}&\dTo_{\De_*I_{abs,w}}\\
&& \De_*(\Om_X^\bullet,\we dw)
\end{diagram}
\end{lem}

\Pf . The proof is similar to that of Lemma \ref{key-lem}. We use the fact that the dualization of the universal Atiyah class 
induces an action of $p_1^*L_w$ on $\OO_\De^{\wt{w}}$ 
so we get a commutative diagram 
\begin{diagram}
p_1^*L_w\ot \OO_\De^{\wt{w}}&\rTo{\id\ot\can}&p_1^*L_w\ot \De_*\De^*\OO_\De^{\wt{w}}\\
\dTo{\hat{\at}^{\univ}}&&\dTo{}\\
\OO_\De^{\wt{w}}&\rTo{\can}&\De_*\De^*\OO_\De^{\wt{w}}
\end{diagram}
and then use the iteration of this action to get a commutative diagram
\begin{diagram}
p_1^*S^n(L_w)\ot \OO_\De^{\wt{w}}&\rTo{\id\ot\can}&p_1^*S^n(L_w)\ot \De_*\De^*\OO_\De^{\wt{w}}\\
\dTo{\exp(\at^{\univ})}&&\dTo{}\\
\OO_\De^{\wt{w}}&\rTo{\can}&\De_*\De^*\OO_\De^{\wt{w}}
\end{diagram}
Finally, composing with $\De_*\vareps:\De_*\De^*\OO_\De^{\wt{w}}\to \OO_\De^{\wt{w}}$
and dualizing we get the result.
\ed

\subsection{Proof of Theorem A}
Let us first consider the case $w=0$. We can view the commutative triangle of Lemma \ref{key-lem} 
as the triangle of Fourier-Mukai functors from $D(X)$ to $D(X)$ (where to a kernel $K$ in $D(X\times X)$
we associate the functor $Rp_{1*}(K\ot p_2^*(\cdot))$. 
Applying these functors to an object $E\in D(X)$, we get a commutative triangle of the form
\begin{equation}\label{canE-triangle}
\begin{diagram}
E &\rTo{\can_E}& \De^*\OO_\De\ot E\\
&\rdTo{\exp(\at(E))}&\dTo_{I_{abs}\ot \id}\\
&& \Om_X^\bullet\ot E
\end{diagram}
\end{equation}
We claim that the morphism $\can_E$ is obtained by applying $\De^*$ to the canonical morphism 
$$\eta_E:p_1^*E\to \De_*E\simeq \OO_\De\ot p_2^*E.$$ 

Indeed, first let us observe that for any $F,G\in D(X\times X)$ we have a commutative triangle
\begin{diagram}
F\ot G &\rTo{a_F\ot \id_G}& \De_*\De^*F\ot G\\
&\rdTo{a_{F\ot G}}&\dTo{\sim}\\
&&\De_*\De^*(F\ot G)
\end{diagram}
where $a_F:F\to \De_*\De^*F$ is the adjunction map, and the vertical arrow is the composition
of the natural isomorphisms
$$\De_*\De^*F\ot G\simeq \De_*(\De^*F\ot \De^*G)\simeq \De_*\De^*(F\ot G).$$ 
Applying this to $F=\OO_\De$ and $G=p_2^*E$ we get the commutativity of the triangle in
the diagram
\begin{diagram}
p_1^*E&\rTo{\eta_E}&\OO_\De\ot p_2^*E&\rTo{\can\ot\id}&\De_*\De^*\OO_\De\ot p_2^*E\\
\dTo{}&&\dTo{}&\ldTo{\sim}\\
\De_*\De^*p_1^*E&\rTo{\De_*\De^*\eta_E}&\De_*\De^*(\OO_\De\ot p_2^*E)
\end{diagram}
Note that here the vertical arrows are the adjunction maps, so the square in the above diagram is also
commutative.
Using the adjointness of $(p_1^*,p_1^*)$, we get a commutative diagram
\begin{diagram}
E&\rTo{\sim}&Rp_{1*}(\OO_\De\ot p_2^*E)&\rTo{Rp_{1*}(\can\ot\id)}&\De^*\OO_\De\ot p_2^*E\\
\dTo{\sim}&&\dTo{}&\ldTo{\sim}\\
\De^*p_1^*E&\rTo{\De^*\eta_E}&\De^*(\OO_\De\ot p_2^*E)
\end{diagram}
By definition, the composition of arrows in the first row is $\can_E$, and our claim follows.

This implies that $\can_E$ corresponds by dualization to the sheafified boundary-bulk map
$$E\ot E^\vee\to \De^*\OO_\De=\und{HH}_*(X)$$
(obtained as $\De^*$ of the evaluation map $E\boxtimes E^\vee\to \OO_\De$).
Composing with $I_{abs}$ and using commutativity of \eqref{canE-triangle}
we get its expression in terms of $\exp(\at(E))$.

Now we can repeat the same argument in the case of matrix factorizations.
For $E\in\MF(X, w)$, using Lemma \ref{key-w-lem}, we get a commutative triangle
\begin{diagram}
E &\rTo{\can_E}& \De^*\OO_\De^{\wt{w}}\ot E\\
&\rdTo{\exp(\at(E))}&\dTo_{I_{abs,w}\ot\id}\\
&& (\Om_X^\bullet,\we dw)\ot E
\end{diagram}
where the morphism $\can_E$ is obtained by applying $\De^*$ to $\eta_E$.
Since $\ev_E$ corresponds to $\eta_E$ by dualization,
this implies that the sheafified boundary-bulk map,
$$\De^*(\ev_E):\und{\End}(E)\to \De^*(\OO_\De^{\wt{w}}),$$
is given by $x\mapsto \str(\exp(\at(E)) \cdot x )$. 
\ed

\section{Graded case}
\label{graded-sec}

\subsection{Basics}

\subsubsection{Category of $\G_m$-equivariant matrix factorizations}

Let $X$ be a smooth scheme equipped with a $\G_m$-action, and let $W$ be a regular function on $X$ satisfying
$$W(\la x)=\la W(x)$$
for $\la\in\G_m$. 

Let $\chi$ denote the identity character of $\G_m$.
The category $\MF_{\G_m}(X,W)$ of $\G_m$-equivariant matrix factorizations of $W$ has as objects
$\G_m$-equivariant $\Z_2$-graded bundles $E=E_0\oplus E_1$ on $X$ equipped with $\OO_X$-linear maps
$$\de_1:E_1\to E_0, \ \ \de_0:E_0\to E_1\ot \chi,$$
such that $\de_0\circ \de_1=\de_1\circ \de_0=W$. In order to define morphisms it is more convenient to replace $(E,\de)$
with the $\Z$-graded bundle equipped with a degree one endomorphism $\de$ such that $\de^2=W$,
$$C(E)=C(E,\de): \ldots E_1\rTo{\de_1} E_0 \rTo{\de_0} E_1\ot \chi \rTo{\de_1} E_0\ot \chi \to\ldots$$
where the $\Z$-grading on $C(E)$ is determined by $C(E)_0=E_0$.
This complex is equipped with a chain {\it $2$-quasi-periodicity} isomorphism
$$\a_E:C(E)\simeq C(E)\ot\chi[-2].$$
Now for a pair of matrix factorizations $E,F$ we define the $\Z$-graded complex of $\G_m$-equivariant bundles,
$\und{\Hom}(E,F)$ as a subcomplex in the sheafified internal $\Hom$-complex of the corresponding $2$-quasi-periodic complexes,
$\und{\Hom}(C(E),C(F))$, consisting of morphisms
respecting isomorphism $\a_E$ and $\a_F$. 

Then the morphisms from $E$ to $F$ are defined as
$$\Hom(E,F):=R\Ga(X,\und{\Hom}(E,F))^{\G_m}.$$
More precisely, here $R\Ga$ should be replaced by some functorial multiplicative resolution.
The resulting category $\MF_{\G_m}(X,W)$ is a $\Z$-graded dg-category, unlike the usual category of matrix factorizations
which is only $\Z_2$-graded (see \cite[Sec.\ 1]{PV-cohFT} for more details).

Note that the complex $\und{\Hom}(E,F)$ is still $2$-quasi-periodic, so in fact we have
$$\und{\Hom}(E,F)=C(\und{\Hom}^{mf}(E,F)),$$
where $\und{\Hom}^{mf}(E,F)$ is a $\G_m$-equivariant matrix factorization of $0$ on $X$.

Note that in the case $W=0$ we can associated with every bounded $\Z$-graded complex of $\G_m$-equivariant vector bundles $(V^\bullet,d)$,
a $\G_m$-equivariant matrix factorization $\mf(V^\bullet)$ of $0$, given by
$$\mf(V^\bullet)_0=\bigoplus_{n\in \Z} V^{2n}\ot \chi^{-n}, \ \ \mf(V^\bullet)_1=\bigoplus_{n\in \Z} V^{2n-1}\ot \chi^{-n}.$$
Note that the corresponding $\Z$-graded complex $C(\mf(V^\bullet))$ is simply
$$C(\mf(V^\bullet))=\bigoplus_{n\in\Z} V\ot \chi^n[-2n].$$

\begin{ex}\label{D(X)-Gm-mf-ex}
In the case when the action of $\G_m$ on $X$ is trivial, one can easily see that the composed functor from 
$D(X)_{\Z}$, the usual $\Z$-graded derived category of coherent sheaves on $X$,
$$D(X)_{\Z}\to D_{\G_m}(X)_\Z\rTo{\mf} \MF_{\G_m}(X,0),$$
where the first functor equips a complex in $D(X)_{\Z}$ with the trivial $\G_m$-action,
is an equivalence (see \cite[Sec.\ 1.2]{PV-cohFT}).
\end{ex}

The above definitions also make sense for more general categories of quasicoherent (resp., coherent)
matrix factorizations.

%For a $\G_m$-equivariant matrix factorization $E$ of $0$ on $X$, we can consider the push-forward
%$R\pi_*E$, where $\pi:X\to pt$, as a $\G_m$-matrix factorization of $0$ on the point.
%We have an identification
%$$C(R\pi_*E)\simeq R\Ga(X,C(E))$$
%as $\G_m$-equivariant sheaves on a point.
%In particular, for $E=\mf(V^\bullet)$, we have

\subsubsection{Tensor product of graded matrix factorizations}

If $W$ and $W'$ are functions satisfying $W(\la x)=\la W(x)$, $W'(\la x)=\la W'(x)$,
then we have a natural operation of tensor product
$$\otimes: \MF_{\G_m}(X,W)\times \MF_{\G_m}(X,W')\to \MF_{\G_m}(X,W+W'),$$
which is uniquely determined by 
\begin{diagram}
C(E\ot F)={\rm equalizer}(C(E)\ot C(F) & \pile{ \rTo^{\a_E\ot \id} \\ \rTo_{\id\ot \a_F}} & C(E)\ot C(F)\ot \chi[-2]).
\end{diagram}
More explicitely,
$$(E\ot F)_0=E_0\ot F_0\oplus E_1\ot F_1\ot\chi, \ \ (E\ot F)_1=E_0\ot F_1\oplus E_1\ot F_0$$
(see \cite[Sec.\ 1.1]{PV-cohFT}).

In the particular case $W=0$ the tensor product gives a structure of a symmetric monoidal category on $\MF_{\G_m}(X,0)$.

There is also a natural duality functor
$$\MF_{\G_m}(X,W)^{op}\to \MF_{\G_m}(X,-W): E\mapsto E^\vee,$$
such that 
$$C(E^\vee)\simeq C(E)^\vee,$$
with $\a_{E^\vee}$ induced by $\a_E$. One has
an isomorphism of $\G_m$-equivariant matrix factorizations of $0$,
$$\und{\Hom}^{mf}(E,F)\simeq F\ot E^\vee$$
(see \cite[Lem.\ 1.1.6]{PV-cohFT}).

Note that we have a natural forgetful functor
\begin{equation}\label{Gm-forget-fun}
\MF_{\G_m}(X,0)\to \MF(X,0)=D(X)_{\Z_2}
\end{equation}
where $D(X)_{\Z_2}$ is the $\Z_2$-graded derived category of coherent sheaves on $X$.
It is easy to see that this functor is compatible with tensor products.

\subsection{Atiyah classes}

Let us consider $\wt{W}=W\ot 1-1\ot W$ as a function on $X^2$. Note that it still satisfies
$\wt{W}(\la x)=\la\wt{W}(x)$, where we equip $X^2$ with the diagonal $\G_m$-action.
We can associate with each object $K$ of $\MF_{\G_m}(X^2,\wt{W})$, with support proper with respect to $p_1$, a Fourier-Mukai type
functor $Rp_{1*}(K\ot p_2^*(\cdot)$ from $\MF_{\G_m}(X,W)$ to itself.

Since $\wt{W}|_{\De(X)}=0$, we have natural functors
\begin{equation}\label{diagonal-Gm-W-functors}
\De_*:\MF_{\G_m}(X,0)\to \MF_{\G_m}(X^2,\wt{W}), \ \ \De^*:\MF_{\G_m}(X^2,\wt{W})\to \MF_{\G_m}(X,0).
\end{equation}
We can view $\OO_X$ sitting in degree $0$ as a $\G_m$-equivariant matrix factorization of $0$ via the functor $\mf$.
We denote by $\OO_\De^{\wt{W}}$ the corresponding object $\De_*(\OO_X)\in \MF_{\G_m}(X^2,\wt{W})$.
Equivalently, 
$$\OO_\De^{\wt{W}}=\mf(\OO_\De).$$

Next, let us define a $\G_m$-equivariant matrix factorization $\OO_{\De,\wt{W}}^{(2)}$ of $\wt{W}$ by
$$C(\OO_{\De,\wt{W}}^{(2)})=[\ldots \OO_{\De^{(2)}}\rTo{-1} \OO_\De \rTo{-dW} \OO_{\De^{(2)}}\ot \chi\to\ldots],$$
where $\OO_{\De^{(2)}}$ sits in degree $0$.

Also, let us define a $\G_m$-equivariant matrix factorization $L^\vee_W$ of $0$, by
$$L^\vee_W=\mf([\OO_X\rTo{dW} \Om_X^1\ot\chi]),$$
where the two-complex is placed in degrees $[0,1]$.
Then we have a natural exact sequence of $\G_m$-equivariant matrix factorizations of $\wt{W}$,
$$0\to \De_*L^\vee_W[-1]\rTo{\phi} \OO_{\De,\wt{W}}^{(2)} \rTo{\psi} \OO_\De^{\wt{W}}\to 0,$$
where the map $\phi$ has as components the identity map on $\OO_\De$ and
the natural embedding $\De_*\Om^1_X\to \OO_{\De^{(2)}}$, while the nontrivial component of $\psi$ is given by
the natural projection $\OO_{\De^{(2)}}\to \OO_\De$.

From the above exact sequence we get a morphism in $\MF_{\G_m}(X^2,\wt{W})$,
$$\hat{\at}^{\univ}: \OO_\De^{\wt{W}}\to \De_*L^\vee_W\simeq p_1^*L^\vee_W\ot \OO_\De^{\wt{W}}.$$
Applying the corresponding Fourier-Mukai functors, we get for every $\G_m$-equivariant matrix factorization $E$ of $W$,
a morphism
$$\hat{\at}(E): E\to L^\vee_W\ot E.$$

\begin{ex}
Assume that the action of $\G_m$ on $X$ is trivial and $W=0$. Then $\MF_{\G_m}(X,0)$ is identified with $D(X)_{\Z}$ 
(see Example \ref{D(X)-Gm-mf-ex}). It is easy to see that in this case $\hat{\at}(E)=\id_E+\at(E)$, where
$\at(E):E\to \Om^1_X[1]\ot E$ is the usual Atiyah class of $E$.
\end{ex}

\subsection{Lie algebra structure, HKR, and the boundary-bulk map}

We are going to equip $L^\vee_W$ with a Lie algebra structure in the symmetric monoidal category $\MF_{\G_m}(X,0)$,
such that the embedding $\mf(\OO_X)\to L^\vee_W$ is in the center.

For this we start by applying the construction of the Atiyah class to $L^\vee_W$ viewed as a matrix factorization of $0$:
$$\hat{\at}_{L^\vee_W}: L^\vee_W\to L^\vee_0\ot L^\vee_W.$$
Next, we observe that there are canonical chain maps of complexes placed in degrees $[0,1]$,
$$[\OO_X\rTo{0}\Om_X^1]\to [0\to \Om_X^1\ot\chi]\to [\OO_X\rTo{dW}\Om_X^1\ot\chi].$$
Applying the functor $\mf$ we get a canonical morphism in $\MF_{\G_m}(X,0)$,
$$\varphi_W:L^\vee_0\to L^\vee_W.$$
Now we define the dual of the Lie bracket as the composition
$$L^\vee_W\rTo{\hat{\at}_{L^\vee_W}} L^\vee_0\ot L^\vee_W\rTo{\varphi_W\ot\id} L^\vee_W\ot L^\vee_W.$$

Next, using adjunction, from $\hat{\at}^{\univ}$, we get a morphism in $\MF_{\G_m}(X,0)$,
$$\iota_W: L_W\to Rp_{1*}\und{\End}^{mf}(\OO_\De^{\wt{W}})=:\und{HH}^*(\MF_{\G_m}(X,W))$$
(recall that $\und{\End}^{mf}(\OO_\De^{\wt{W}})$ is a $\G_m$-equivariant matrix factorization of $0$ on $X^2$).
Using the product in $\und{HH}^*(\MF_{\G_m}(X,W))$ we obtain a map of $\G_m$-equivariant matrix factorizations of $0$
$$I^{abs,\G_m}:S^n(L_W)\to \und{HH}^*(\MF_{\G_m}(X,W)),$$
where $n=\dim X$.

On the other hand, setting
$$\und{HH}_*(\MF_{\G_m}(X,W)):=\De^*\OO_\De^{\wt{W}}\in \MF_{\G_m}(X,0),$$
we have a natural action of $\und{HH}^*(\MF_{\G_m}(X,W))$ on $\und{HH}^*(\MF_{\G_m}(X,W))$
(understood in terms of the tensor structure on $\MF_{\G_m}(X,0)$), and a canonical functional
$$\vareps:\und{HH}_*(\MF_{\G_m}(X,W))\to \mf(\OO_X),$$
coming from the adjunction for the functors \eqref{diagonal-Gm-W-functors}.
Thus, iterating the action of $L_W$ on $\und{HH}_*(\MF_{\G_m}(X,W))$ and applying $\vareps$, we get similarly to the $\Z_2$-graded case
a map
$$I_{abs,\G_m}:\und{HH}_*(\MF_{\G_m}(X,W))\to S^nL^\vee_W.$$

Since the maps $I^{abs,\G_m}$ and $I_{abs,\G_m}$ lift the previously defined maps $I^{abs}$ and $I_{abs}$ using the
forgetful functor \eqref{Gm-forget-fun}, they are quasi-isomorphisms.

Note that we have an isomorphism defined in the same way as \eqref{Sn-isom},
$$S^nL^\vee_W\simeq \mf(S^n[\OO_X\rTo{dW} \Om^1_X\ot\chi])\simeq \mf({\bigwedge}^\bullet(\Om^1_X\ot \chi),\we dW).$$
Thus, iterating $\hat{\at}(E)$, as in the $\Z_2$-graded case, we get for $E\in\MF_{\G_m}(X,W)$, a map
$$\exp(\at(E)): E\to \mf({\bigwedge}^\bullet(\Om^1_X\ot \chi),\we dW)\ot E.$$

Now all the previous arguments generalize immediately to the graded case and give the following analog of Theorem A.

\medskip

\noindent
{\bf Theorem B}. {\it Assume that $W=0$ on the critical locus of $W$ (set-theoretically). 
There is a natural identification 
\begin{equation}\label{Gm-HKR-eq}
HH_*(\MF_{\G_m}(X,W))\simeq H^*(X,{\bigwedge}^\bullet(\Om^1_X\ot \chi),\we dW)^{\G_m}.
\end{equation}
%\footnote{BK: Is $R\Gamma (X, (\Om_X^\bullet \we dw ))$ instead of $(\Om_X^\bullet \we dw )$?}
Under this identification, the categorical boundary-bulk map for $E\in \MF_{\G_m}(X,W)$, 
%\footnote{BK: in $D(k)$?}
$$\Hom^*_{\MF_{\G_m}(X,W)}(E,E)\to HH_*(\MF_{\G_m}(X,W)),$$
is equal to the map induced on the $\G_m$-invariant part of hypercohomology by the map
$$\und{\Hom}^{mf}(E,E)\to \mf({\bigwedge}^\bullet(\Om^1_X\ot \chi),\we dW): x\mapsto \str(\exp(\at(E)) \cdot x)$$
in $\MF_{\G_m}(X,0)$.
}

\end{document}